\documentclass[10pt]{article}

\input xy
\xyoption{all}
\usepackage{amsmath}
\usepackage{amsfonts}
\usepackage{amssymb}
\usepackage{stmaryrd}
\usepackage{graphicx,epsfig}
\usepackage{bbm,color}
\usepackage{dbnsymb}
\usepackage{amsthm}

\newtheorem{thm}{Theorem}[section]
\newtheorem{defn}[thm]{Definition}
\newtheorem{lem}[thm]{Lemma}
\newtheorem{quest}[thm]{Question}
\newtheorem{prop}[thm]{Proposition}
\newtheorem{rem}[thm]{Remark}
\newtheorem*{prop*}{Proposition}

\renewcommand{\phi}{\varphi}
\newcommand{\sg}{\langle}
\newcommand{\sd}{\rangle}
\newcommand{\boA}{\mathcal{A}}
\newcommand{\boB}{\mathcal{B}}
\newcommand{\boD}{\mathcal{D}}
\newcommand{\boH}{\mathcal{H}}
\newcommand{\boM}{\mathcal{M}}

\newcommand{\boC}{\mathcal{C}}
\newcommand{\boI}{\mathcal{I}}
\newcommand{\boG}{\mathcal{G}}
\newcommand{\boL}{\mathcal{L}}

\newcommand{\Q}{\mathbb{Q}}
\newcommand{\Z}{\mathbb{Z}}

\newcommand{\alg}[3]{\vphantom{#2}_{#1} #2_{#3}}
\newcommand{\mat}[1]{\begin{pmatrix} #1 \end{pmatrix}}

\newcommand{\dd}{\text{d}}
\newcommand{\arete}{\!\!\frown}
\newcommand{\Zrat}{Z^{\rm rat}}
\newcommand{\irat}{\int^{\rm rat}}

\newcommand{\ba}{\overline}
\DeclareMathOperator{\tp}{ST}
\DeclareMathOperator{\id}{Id}
\DeclareMathOperator{\lift}{Lift}
\DeclareMathOperator{\hair}{Hair}
\DeclareMathOperator{\Tr}{Tr}

\begin{document}
\title{Surgery on a single clasper and the 2-loop part of the Kontsevich integral
\footnote{\emph{keywords}: finite type invariants, Kontsevich integral, rationality, S-equivalence. 
\emph{2000 Mathematics Subject Classification:} 57M27.}}
\author{Julien March\'e
\footnote{Institut de Math\'ematiques de Jussieu, \'Equipe ``Topologie et G\'eom\'etries Alg\'ebriques''
Case 7012, Universit\'e Paris VII, 75251 Paris CEDEX 05, France. 
\quad e-mail: \texttt{marche@math.jussieu.fr}}}
\date{October 2004}
\maketitle

\begin{abstract}
We study the 2-loop part of the rational Kontsevich integral of a knot in an integer homology sphere. We give a general formula which explains how the 2-loop part of the Kontsevich integral of a knot changes after surgery on a single clasper whose leaves are not linked to the knot. As an application, we relate this formula with a conjecture of L. Rozansky about integrality of the 2-loop polynomial of a knot.

\end{abstract}

\section{Introduction}

The Kontsevich integral of a knot in $S^3$ was introduced in 1993 in \cite{kont}. It is a an invariant as powerful as it is difficult to compute. It is organized in a series of diagrams whose $n$-th term $z_n$ contains all diagrams with loop degree $n$. The loop degree 1 of the Kontsevich integral is determined by the Alexander polynomial, hence the difficulty of the computation begins with the 2-loop part of the Kontevich integral. This invariant may be seen as an equivariant Casson invariant because it contains the Casson invariants of the ramified coverings of $S^3$ along a knot $K$ as it was shown in \cite{lift}.


Let $z_2(M,K)$ be the 2-loop part of the Kontsevich integral of a pair $(M,K)$ where $K$ is a knot in an integer homology sphere $M$. It is an element of $\boA^{\Delta}_{2}$, i.e. it is presented as a sum of graphs $\ThetaGraph$ and $\dumbbell$ colored by rational expressions such as $\frac{P}{\Delta}$ where $P\in\Q[t,t^{-1}]$ and $\Delta$ is the Alexander polynomial of $(M,K)$. We refer to \cite{rat} for a precise construction of this invariant. We note $\boI^{\Delta}_2$ the subgroup of $\boA^{\Delta}_{2}$ spanned by graphs such that all polynomials $P$ have integer coefficients (see section 2 for precise definitions).

We give a formula for $z_2((M,K)^G)-z_2(M,K)$ where $(M,K)^G$ is the surgery of $(M,K)$ on a single clasper $G$ in $M\setminus K$ with leaves homologous to 0 (see \cite{hab}, \cite{clasp}, or \cite{gr}). This formula uses a kind of equivariant Milnor number of the clasper. It gives a new tool for the computation of 2-loop parts and is very convenient for studying the relation between 2-loop part and $S$-equivalence.

We will deduce from this formula that if $(M,K)$ and $(M',K')$ are $S$-equivalent, then $z_2(M,K)-z_2(M',K')$ lives in $\frac{1}{2}\boI^{\Delta}_2$.

In \cite{roz},  L.Rozansky conjectured that if $K$ is a knot in an integer homology sphere $M$, then $z_2(M,K)$ belongs to $\frac{1}{2}\boI^{\Delta}_2$. We give a formula for $z_2(M,K)\mod \frac{1}{2}\boI^{\Delta}_2$ from the $S$-equivalence class which is conjectured to give 0. We use it to prove that $z_2(M,K)\in \frac{1}{12}\boI^{\Delta}_2$.

{\bf Aknowlegments:} We wish to thank T. Ohtsuki, A. Kricker, P. Vogel and Y. Tsutsumi for useful comments.

\section{General facts about the 2-loop part of the Kontsevich integral}

Let $M$ be a rational homology sphere and $K$ be a knot homologous to 0 in $M$. Then, the rational Kontsevich integral constructed in \cite{rat} is an invariant of the pair $(M,K)$ lying in a space $\boA^{\Delta}$. 

This space is made of closed trivalent diagrams whose edges are colored by elements such as $\frac{P}{\Delta}$ where $P\in\Q[t,t^{-1}]$ and $\Delta=\Delta(M,K)$ is the Alexander polynomial of the pair $(M,K)$. In order to define this space, we use a more general construction of space of diagrams.

\begin{defn}
Let $k$ be a ring, $H$ a co-commutative Hopf algebra over $k$ (or a Hopf algebra up to completion) and $E$ a $H$-module with an involution noted $x\mapsto \ba{x}$ and a ``unit'' noted $1_E$.
We define the space $\boD(E)$ of diagrams colored by $E$ by generators and relations:

\begin{itemize}
\item
generators are couples $(\Gamma,f)$ where $\Gamma$ is a finite trivalent diagram with a cyclic ordering of the edges incoming at each vertex and $f$ is a collection of pairs $(\omega_e,x_e)_{e}$ for all edges of $\Gamma$ where $\omega_e$ is an orientation of the edge and $x_e$ is an element of $E$.

\item
the relations are the following:
\begin{enumerate}
\item
changing orientation of an edge - \\
let $f'$ be the family $f$ except for some edge $e$ where we have replaced $(\omega_e,x_e)$ by $(-\omega_e,\ba{x_e})$. Then we identify $(\Gamma,f)$ and $(\Gamma,f')$.
\item
sliding an element of $H$ - \\
At each vertex of $\Gamma$, we introduce an action of $H\otimes H\otimes H$ by multiplying the corresponding edges by scalars in $H$ in a way compatible with orientation.
If $h\in H$ and $\Delta h =\sum_{h',h''}h'\otimes h''$, we identify the action of the elements of figure \ref{sliding}.
\begin{figure}[htbp]
\begin{center}
\input{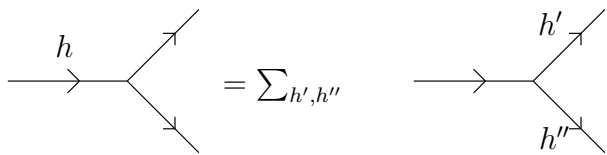}
\caption{Sliding relation}
\label{sliding}
\end{center}
\end{figure}
\item
global symetry - \\
if there is a isomorphism $(\Gamma,f)\simeq (\Gamma',f')$ which respects local orderings then, we identify the resulting generators.
\item
antisymetry - \\
Let $\Gamma'$ be the graph $\Gamma$ with the reversed ordering at one vertex. Then we set $(\Gamma',f)=-(\Gamma,f)$.
\item 
relation IHX - \\
If some edge is colored by $1_E$, we may perform a $\IGraph=\HGraph-\XGraph$ relation at this edge.
\end{enumerate}
\end{itemize}
\end{defn}

The main example is the following: let us take $k=\Q$, $H=\Q[t,t^{-1}]$ and $E_{\Delta}=\{\frac{P}{\Delta},P\in\Q[t,t^{-1}]\}$ where the involution maps $t$ to $t^{-1}$. We obtain a space noted $\boA^{\Delta}$ where $\Zrat(M,K)$ lives.

If we take $H=E=\Q[[h]]$, we obtain the usual space of open diagrams which we note $\boB$. The substitution $t\mapsto \exp(h)$ gives a map from $E_{\Delta}$ to $\Q[[h]]$ which induces a map from $\boA^{\Delta}\to \boB$  which is noted $\hair$.

Let us give a complete proof of the following well-known fact:
\begin{prop}\label{deuxloop}
  The space $\boA^{\Delta}_2$ of 2-loop rational diagrams with denominators $\Delta$ is isomorphic to $\Q[H^1(\ThetaGraph)]_{Aut(\ThetaGraph)}$ as a $\Q$-vector space. Moreover, the map $\hair:\boA^{\Delta}_2\to \boB$ is injective.
\end{prop}
\begin{proof}

For any trivalent diagram $D$ with orderings at vertices, we note  $\boA^{\Delta}(D)$ the space of its colorings by rational functions with denominator $\Delta$, that is a tensor product of copies of $E$ for all edges, modulo the sliding relations. The automorphism group of $D$ acts on $\boA^{\Delta}(D)$ by changing the sign any time some orientation is changed at a vertex. We note $\boA^{\Delta}(D)^{Aut(D)}$ the space of coinvariants under this action. It is clear that $\boA^{\Delta}_2=\boA^{\Delta}(\ThetaGraph)^{Aut(\ThetaGraph)}\oplus\boA^{\Delta}(\dumbbell)^{Aut(\dumbbell)} /(IHX)$ because there are only 2 connected diagrams of loop degree 2, namely $\dumbbell$ and $\ThetaGraph$.

Let $i$ be the following map induced by the inclusion:
$$\xymatrix{
\boA^{\Delta}(\ThetaGraph)^{Aut(\ThetaGraph)} \ar[r] &\boA^{\Delta}(\ThetaGraph)^{Aut(\ThetaGraph)}\oplus\boA^{\Delta}(\dumbbell)^{Aut(\dumbbell)} /(IHX)
}.$$

We want to show that it is an isomorphism. Hence, we define a retraction $\phi$ which is the identity on the component $\boA^{\Delta}(\ThetaGraph)^{Aut(\ThetaGraph)}$ and which is defined on $\boA^{\Delta}(\dumbbell)^{Aut(\dumbbell)}$ in the following way.

Let $D$ be the diagram $\frac{P}{\Delta}\overset{\frac{R}{\Delta}}{\bigcirc{\!-\!}\bigcirc}\frac{Q}{\Delta} \in \boA^{\Delta}(\dumbbell)$. We can write $\frac{R}{\Delta}=R(1)+\frac{R-R(1)\Delta}{\Delta}=R(1)+\frac{(t-1)R'}{\Delta}$ for some $R'\in\Q[t,t^{-1}]$.

By applying a sliding relation, we see that a diagram $\dumbbell$ containing $t-1$ on its central edge vanishes in $\boA^{\Delta}(\dumbbell)$.

Then, we can suppose that this edge is colored by 1 and we may write an (IHX) relation. This relation identifies the diagram $D$ with an element of $\boA^{\Delta}(\ThetaGraph)$, which gives the following relation in 
$\boA^{\Delta}_2$.

\begin{equation}\label{application}
\frac{P}{\Delta}\bigcirc{\!\!-\!\!}\bigcirc\frac{Q}{\Delta}=\theta(\frac{P}{\Delta},\frac{\ba{Q}-Q}{\Delta}).
\end{equation}

In this equation, $\theta(\frac{P}{\Delta},\frac{Q}{\Delta},\frac{R}{\Delta})$ is a graph $\ThetaGraph$ whose three edges are oriented in the same direction and are colored by $\frac{P}{\Delta},\frac{Q}{\Delta}$ and $\frac{R}{\Delta}$. If some function is missing, it means that it may be replaced by 1.

If we suppose that some edge of $\ThetaGraph$ has no colors and that we apply to it an (IHX) relation, we get back the relation \eqref{application}. Hence, the only (IHX) relation is the relation \eqref{application} which identifies any diagram of $\boA^{\Delta}(\dumbbell)$ with a diagram of $\boA^{\Delta}(\ThetaGraph)$.

Then, we can set $\phi(\frac{P}{\Delta}\bigcirc{\!\!-\!\!}\bigcirc\frac{Q}{\Delta})=\theta(\frac{P}{\Delta},\frac{\ba{Q}-Q}{\Delta})$. This map respects the (IHX) relations by definition. We check easily that it induces a map on coinvariants. The existence of this retraction prooves that $i$ is an isomorphism. 
We can replace $\boA^{\Delta}_2$ by $\boA^{\Delta}(\ThetaGraph)^{Aut(\ThetaGraph)}$.

Now, we identify $\boA^{\Delta}(\ThetaGraph)^{Aut(\ThetaGraph)}$ et $\Q[H^1(\ThetaGraph)]_{Aut(\ThetaGraph)}$ on the following way:

Suppose that $x=\sum_i \frac{P_i}{\Delta}\otimes \frac{Q_i}{\Delta}\otimes \frac{R_i}{\Delta} \in \boA^{\Delta}(\ThetaGraph)$. We identify $t\otimes 1\otimes 1$, $1\otimes t\otimes 1$ and $1\otimes 1\otimes t$ with the 3 classes of $H^1(\ThetaGraph,\Z)$ corresponding to the edges.

We write $\ba{x}=\sum_i\sum\limits_{\tau\in \Z/2\times S_3} \tau. (P_i\otimes Q_i\otimes R_i)\in\Q[H^1(\ThetaGraph,\Z)]_{Aut(\ThetaGraph)}$. This map is clearly an isomorphism.

Let us proove now the injectivity of the $\hair$ map on $\boA^{\Delta}_{2}$.

The $\hair$ map corresponds to the substitution $t\mapsto \exp(h)$. It is sufficient to show that the $\hair$ map is injective on $\boA^{\Delta}(\ThetaGraph)^{Aut(\ThetaGraph)}$. By symetrization and multiplication by $\Delta^{\otimes 3}$, it is sufficient to show that the map $\hair:\Q[H^1(\ThetaGraph,\Z)]\to\Q[[H^1(\ThetaGraph,\Z)]]$ is injective. 
But $H^1(\ThetaGraph,\Z)$ is a free module with 2 generators $a$ and $b$, and as $\exp(a)$ and $\exp(b)$ are algebraically independant on $\Q[[a,b]]$, the above map has to be injective.

\end{proof}

Every monomial in $\Q[H^1(\ThetaGraph)]$ is equivalent modulo the action of the group $Aut(\ThetaGraph)$ to a monomial such as $\theta(t^m,t^n)$ for $m$ and $n$ with $0\le 2n\le m$. Moreover, two such monomials are not equivalent. Then we can write in a unique way the polynomial $\ba{z_2}(M,K)$ as a combination of the symetrized monomials $\ba{\theta(t^m,t^n)}$, which gives the most concise formulas. This definition is very closed to the definition of L. Rozansky in \cite{roz}, although there is a subtle difference: in his paper, L. Rozansky considered coefficients of $a^mb^n$ ($0\le 2n\le m$) in $\ba{z_2}(S^3,K)\in\Q[a^{\pm1},b^{\pm 1}]$, which does not give exactly the same basis as the $\ba{\theta(t^m,t^n)}$s.

Let us note $s^mt^n=\ba{\theta(t^m,t^n)}$ and call $m$ its degree. L. Rozansky has proposed the following conjectures for $M=S^3$.

\begin{enumerate}\item
Let $g(K)$ be the genus of $K$ (that is the minimal genus of a Seifert surface of $K$). Then the inequality $\deg \ba{z_2}(S^3,K) \le 2 g(K)$ holds.

\item
The expression $\ba{z_2}(S^3,K)$ is a linear combination of the elements $s^mt^n$ with half integer coefficients for $n=0$ or $m=2n$, and integer coefficients for $0<2n<m$.
\end{enumerate}

To translate this conjecture, we introduce the following definition:
\begin{defn}
Let $k=\Z$, $H=\Z[t,t^{-1}]$ and $E=\{\frac{P}{\Delta},P\in \Z[t,t^{-1}]\}$. We define the space $\boI^{\Delta}$ as the image of $\boD(E)$ in $\boA^{\Delta}$. We will call integer diagrams the elements of $\boI^{\Delta}$.

\end{defn}

The Rozansky conjecture is very closed to the following modified conjecture: for any knot $K$ in an integer homology sphere $M$, we have $z_2(M,K)\in \frac{1}{2}\boI^{\Delta}_2$. We will prove the weaker statement $z_2(M,K)\in \frac{1}{12}\boI^{\Delta}_2$, and we will give some ideas to prove or disprove Rozansky conjecture.

\section{A surgery formula}

\subsection{Notations}
Let $(M,K)$ be a pair formed of a knot $K$ homologous to 0 in a rational homology sphere and $\Gamma$ a trivalent banded graph in the complement of $K$.

We recall that there is a rational integral of $\Gamma$ in $(M,K)$ which is noted $\Zrat_{\alpha}((M,K),\Gamma)$ and which is defined on the following way:

Let $L$ be a link presenting $M\setminus K$ by surgery in the solid torus $\tp$ and $\Gamma'$ be the image of $\Gamma$ in $\tp\setminus L$.

Then, we define
$$\Zrat_{\alpha}((M,K),\Gamma)=\frac{\irat \sigma Z_{\alpha}(\hat{L}\cup \Gamma',\tp)\dd L}{(\irat \sigma Z_{\alpha}(\hat{U}^1,\tp)\dd U^1)^{\sigma^+(L)}
(\irat \sigma Z_{\alpha}(\hat{U}^{-1},\tp)\dd U^{-1})^{\sigma^-(L)}}.$$
In this formula we use the following notations (see \cite{rat} for more details):
\begin{itemize}
\item
$Z(\cdot,\tp)$ is the Kontsevich integral in the solid torus.
\item
$\sigma$ is the inverse Poincar\'e-Birkhoff-Witt map. It is defined for sliced components in a solid torus, that is, components cut along some meridional disc of the solid torus.
\item
$U^n$ is the trivial knot with framing $n$.
\item
$\hat{L}$ means that we have added a normalisation $\nu=Z(U)$ to the corresponding components, where $Z$ is the standard Kontsevich integral.
\item
$\alpha$ is a global normalisation factor: if $\alpha=1$, we obtain the unwheeeled integral, if $\alpha=\nu=Z(U)$, we obtain the standard integral. 
\item
$\irat$ is a formal gaussian integration.
\item
$(\sigma^+(L),\sigma^-(L))$ is the signature of the linking matrix of $L$ in $S^3$.
\end{itemize}

The Fubini theorem implies that if $\Gamma$ is a link in $M\setminus K$ homologous to 0 whose linking matrix is invertible over $\Q$, then the series $\Zrat_{\alpha}((M,K)^{\Gamma})$ is equal to the integral  $\irat \Zrat_{\alpha}((M,K),\hat{\Gamma})\dd \Gamma/N$ where $N$ is the usual normalization depending on the linking matrix of $\Gamma$.

In what follows, we will use this fact to give a formula which describes the action of a surgery over a single clasper on the 2-loop part of the rational Kontsevich integral.

\subsection{The equivariant Milnor number $\mu$}

\begin{figure}[htbp]
\begin{center}
\input{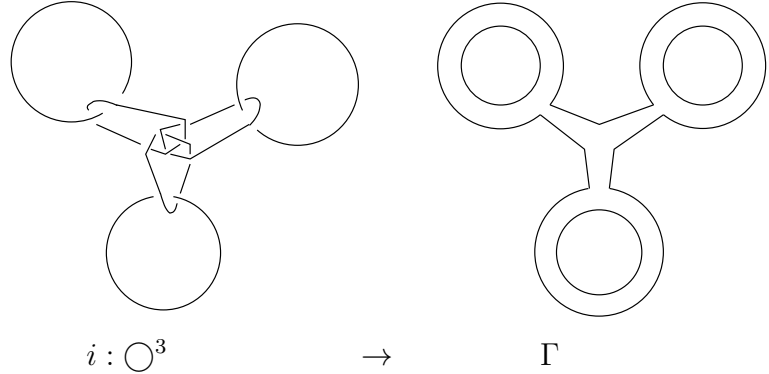}
\caption{Clasper and surgery on a clasper}
\label{clasper}
\end{center}
\end{figure}

Let $\Gamma$ be the banded trivalent graph of the figure \ref{clasper}. Suppose we have an embedding homologous to 0 of this graph in $M\setminus K$.

Let $\boA^{\Delta}(\Gamma)$ be the space of trivalent diagrams colored by $E_{\Delta}$ and lying on the graph $\Gamma$ modulo the standard relations.
The series $\Zrat_{\alpha}((M,K),\Gamma)$ is well defined in $\boA^{\Delta}(\Gamma)$. Let us consider the map $\phi:\boA^{\Delta}(|||)\to\boA^{\Delta}(\Gamma)$ induced by the inclusion of the leaves. It is clearly surjective.

\begin{figure}[htbp]
\begin{center}
\input{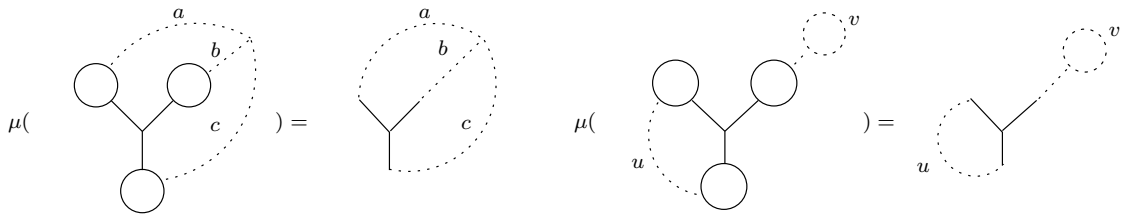}
\caption{The map $\mu$}
\label{mu}
\end{center}
\end{figure}

We define $\mu(\Gamma)$ in the following way:
let us choose a group-like element $X$ in $\boA^{\Delta}(|||)$ such that $\phi(X)=\Zrat_{\alpha}((M,K),\Gamma)$, and then apply the symetrization map $\sigma$ to its 3 components. Consider the diagrams in $X$ of degree 2, which have only 1 trivalent vertex, and one leg for each leaf. By collapsing each leaf, we obtain an element of degree 2 which we note $\mu(\Gamma)$ (see figure \ref{mu})

\begin{lem}
  The map $\mu:\boA^{\Delta}(\Gamma)\to\boA_2^{\Delta}$ is well defined. We may think of $\mu(\Gamma)$ as a contracted equivariant Milnor number of the leaves of $\Gamma$.
\end{lem}

\begin{proof}
It is clear that all diagrams of $\boA^{\Delta}(\Gamma)$ which contribute to $\mu(\Gamma)$ are presented in figure \ref{mu}. 

The ambiguity comes from the choice of the preimage by $\phi$. Two such choices are linked by the sliding of a dashed trivalent vertex through some leaves of $\Gamma$. After the collapsing, we interpret this move as a sliding in a 2-loop rational diagram. This shows that the map $\mu$ is well defined as a map from  $\boA^{\Delta}(\Gamma)$ to $\boA_2^{\Delta}$.
\end{proof}

\subsection{The formula}

We are now able to state the surgery formula.

\begin{prop}\label{claspsimp}
Let $(M,K)$ be a pair formed of a knot homologous to 0 in a rational homology sphere and $\Gamma$ a clasper homologous to 0 in $M\setminus K$. Let $\Gamma'$ be a clasper parallel to $\Gamma$ (pushed of $\Gamma$ thanks to the banding).

We have :
$$z_2((M,K)^{\Gamma})-z_2(M,K)=\frac{1}{2}\sg (M,K),\Gamma\cup\Gamma'\sd + \mu(\Gamma).$$

In this formula, $\sg (M,K),\Gamma\cup\Gamma'\sd$ is the sum of all ways of contracting the leaves of $\Gamma\cup\Gamma'$ using their negative Blanchfield pairing (see \cite{rat}).
\end{prop}

\begin{proof}
Let $L$ be a link in the solid torus $\tp$ which represent the pair $(M,K)$.

Let $\boG=\{F_1,F_2,F_3,A_1,A_2,A_3\}$ be the link with 6 components presented by $\Gamma$ in $\tp\setminus L$ (see figure \ref{milnor}).
The component $F_i$ are the leaves and the components $A_i$ are edges. We call $B_i=F_i\cup A_i$ the arms of the clasper $\Gamma$. Our aim is to determine $z_2([(M,K),G])$. 

\begin{figure}[htbp]
\begin{center}
\input{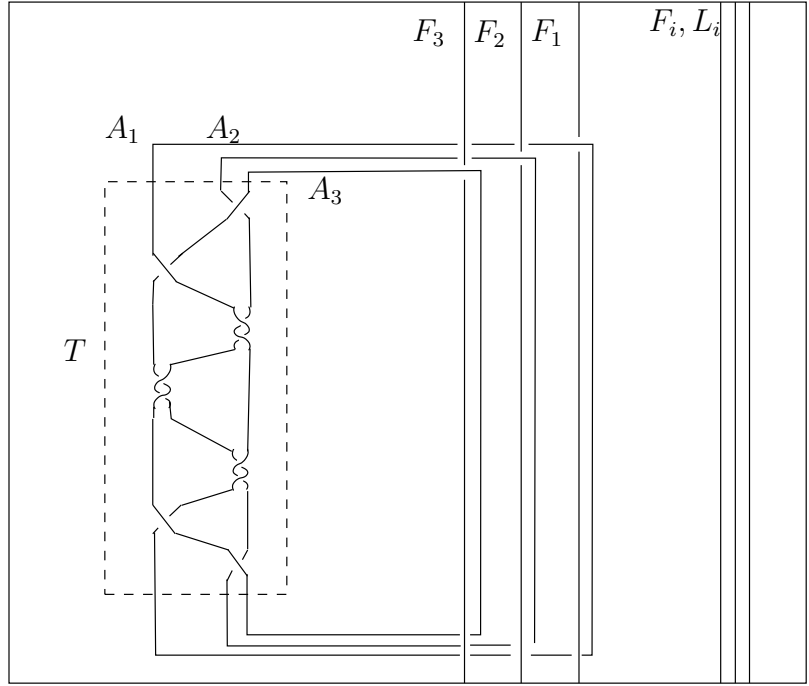}
\caption{Surgery link}
\label{milnor}
\end{center}
\end{figure}

If we perform a surgery on a strict part of $B_1,B_2,B_3$, we do not modify the couple $(M,K)$. This is a fundamental property of borromean surgery. We can write 
$$[(M,K),B_1\cup B_2\cup B_3]=\sum\limits_{I\subset\{1,2,3\}}(-1)^I(M,K)^{B_I}=(M,K)^G-(M,K).$$

We deduce the following formula:
 $$\Zrat((M,K)^G)-\Zrat(M,K)=\sum\limits_{I\subset\{1,2,3\}}(-1)^I\irat \sigma Z(B_I\cup L,\tp)\dd L\dd B_I.$$

We can use this fact by ignoring in $\irat \sigma Z(B_1\cup B_2 \cup B_3\cup L,\tp)\dd L$ all diagrams which are not linked to $B_1,B_2$ and $B_3$.
 The diagrams which produce a 2-loop term are enumerated in the figure \ref{white7}.

\begin{figure}[htbp]
\begin{center}
\input{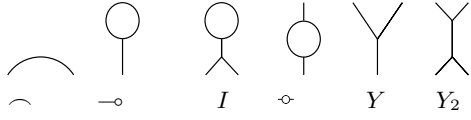}
\caption{Diagrammes producing 2-loop terms}
\label{white7}
\end{center}
\end{figure}

By the preceding remark, we can ignore diagrams such as $\twowheel$ and $I$. We then have to care only on the diagrams $Y$, $Y_2$ and $\multimap$.

We compute explicitely these contributions from the surgery presentation of figure \ref{milnor}. This computation is not obvious because we need an explicit formula for the integral $Z(T)$.

Let us explain this computation. Let $a \in\boA(|||)$ be the diagram made of a single edge lying on the 2 first strands with compatible orientations and $b$ the same diagram lying on the two last strands. We compute $Z(T)$ modulo connected diagrams of degree 4 and more. As the associator $\Phi$ is supposed to be even, we have $\Phi\equiv\exp([a,b]/24)$.

We use Maple to compute $Z(T)$ in the algebra of formal power series in the non commutative variables $a$ and $b$.

\begin{eqnarray*}
Z(T)&=&\exp(-b/2)\exp([a,b]/24)\exp(a/2)\exp(-[a,b]/24)\exp(b)\exp([a,b]/24)\exp(-a)\\
&&\exp(-[a,b]/24)\exp(-b)\exp([a,b]/24)\exp(a/2)\exp(-[a,b]/24)\exp(b/2)\\
&=&\exp([a,b])
\end{eqnarray*}

We have chosen such a $T$ to obtain this simple formula. We could have presented $T$ in a simpler form but then $Z(T)$ would have contained trees with two internal vertices which is not convenient for our purposes.

We deduce the following expression:

\begin{equation}\label{integclasp}
\irat \sigma Z(\boG\cup L,\tp)\dd L = \exp( \frac{1}{2}\sum_{i,j} \overset{w_{ij}}{\alg{F_i}{\arete}{F_j}}+ \sum_i \alg{A_i}{\arete}{F_i}+\multimap (F_i) + \YGraph(F_i)+\YGraph_{A_1,A_2,A_3}+Y_2+\cdots).
\end{equation}
In this formula, $Y_2$ is made of diagrams having the shape of $Y_2$. 
From the computation of $Z(T)$, we know that there are no such diagrams completely colored by the $A_i$s.

The linking matrix of $\boG$ is $M=\mat{W&\id \\ \id & 0}$ where $W$ is the equivariant linking matrix of the leaves of $G$. We have $M^{-1}=\mat{0&\id \\ \id & -W}$.

In the rational integration, we glue the legs $A_i$ themselves with the opposite of the linking matrix, and we glue directly the $A_i$s to the corresponding $F_i$s.

The diagrams $Y_2$ are glued themselves during the integration process. Hence, if some of them contains a $F_i$, it must contain the corresponding $A_i$. But this diagram must also contain the three colors. The only possiblity comes from the diagram colored by $F_i,A_i,A_j,A_k$.
This contribution comes from the multiplication on $F_i$ of diagrams from $Z(T)$ and $\exp_{\#}(\alg{A_i}{\arete}{F_i})$. It creates a diagram $Y_2$ in which the legs $A_i$ and $F_i$ are adjacent. By integration, these two legs are linked, and the diagram vanishes. This proves that we can ignore all $Y_2$ diagrams.

If we perform the integration in the formula \ref{integclasp}, we obtain exactly the result announced in the proposition.

\end{proof}

\section{Integrality}\label{integralite}


\subsection{Integrality properties of $\mu$}

In this part, we will show that $\mu(\Gamma)$ is half-integer before showing that for all pairs $(M,K)$ we have $z_2(M,K)\in \frac{1}{12}\boI^{\Delta}_2$. In both cases, we use a simplified surgery presentation which is described in the following lemma:

\begin{lem}\label{projection}
\begin{enumerate}
\item
Let $\Sigma$ be an annulus such that $\tp=\Sigma\times [0,1]$. Each pair $(M,K)$ made of a link homologous to 0 in a rational homology sphere is presented by surgery on a link in $\tp$ whose projection on $\Sigma$ is the following: each component is trivial with some framing and 2 components cross themselves in at most two points.

\item
Let $\boC$ be the set of sliced links obtained by stacking and juxtaposing the elementary tangles of the figure \ref{classec} up to some symetry. The sign of the linking of the components of $C$ is arbitrary.

\begin{figure}[htbp]
\begin{center}
\input{classec.pstex_t}
\caption{Elements generating the class $\boC$}
\label{classec}
\end{center}
\end{figure}
Each pair $(M,K)$ formed of a knot homologous to 0 in a rational homology sphere is $S$-equivalent to a pair presented by surgery on a link belonging to the class $C$, with some twists added to the components.
\end{enumerate}
\end{lem}

\begin{proof}

Let us start from a surgery presentation of $(M,K)$ by any link $C_i$. We can change any crossing between $C_i$ and $C_j$ by adding a surgery component. This component bounds a disc $D$ which cuts the link in only two points $A$ and $B$ linked in $D$ with a segment $[A,B]$. The whole disc retracts on a neighborhood of this segment. We call special such a component. 
By adding special components, we can trivialize the $C_i$ in such a way that they project on $\Sigma$ on disjoint discs. We can shrink the special components around their segment. Up to isotopy, we can suppose that they project generically on $\Sigma$.

The following step is to eliminate self-crossings of special components. To do this, we use the fact that we can split a special component in a Hopf link of special components. Hence, we can replace each special component by a chain such that all intersections involve different components. We can also suppose that two components cross in two points at most. The move we use is presented in the figure \ref{chaine}

\begin{figure}[htbp]
\begin{center}
\input{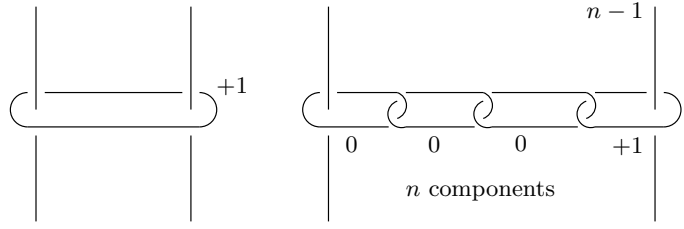}
\caption{Chain of surgery annuli}
\label{chaine}
\end{center}
\end{figure}

This finishes the first part of the lemma.

To show the second part, we remark that the $S$-equivalence class of the pair obtained by surgery on a link $L$ is characterized by the equivariant linking matrix of $L$. This means that if we make a crossing change between a special disc and a component $C_i$, we do not change the $S$-equivalence class. Then, we can suppose that in the neighborhood of the $C_i$s, the tangle looks as the tangle $C$ of figure \ref{classec}.

The link formed of the trivialized components $C_i$ and the special components is convenient and proves the second part of the lemma.
\end{proof}

This lemma allows us to proove a first integrality result concerning the map $\mu$.
\begin{lem}\label{demientier}
  If $K$ is a knot in an integer homology sphere $M$, for all embedding of $\Gamma$ in $M\setminus K$ homologous to 0, $\mu(\Gamma)$ belongs to $\frac{1}{2}\boI^{\Delta}_2$.

\end{lem}

\begin{proof}
Let $\Gamma$ be a clasper in $M\setminus K$ where $M$ is an integer homology sphere.
We can present $(M,K)$ by surgery on a link $\boL$ in $\tp$. Let $F_1,F_2,F_3$ the image of the leaves of $\Gamma$ in $\tp\setminus\boL$. Using the lemma, we can untie $\boL\cup\{F_1,F_2,F_3\}$ using special components in such a way that each component is trivial (and hence bounds a disc $D_i$) and that two components cross themselves in at most two points.

We decompose the link $\boL\cup\{F_1,F_2,F_3\}$ in elementary tangles and we choose a marked point for each component.

It is sufficient to see that all diagrams with shape  $\YGraph$ and $\multimap$ which appear in $\sigma Z(\boL\cup F_1\cup F_2 \cup F_3,\tp)$ are half-integers. These diagrams have degree at most 2. Then, we quotient the space of diagrams by diagrams with degree at least 3 and diagrams with connected components having two trivalent vertices or more.

We distinguish two kinds of diagrams in $Z(\boL\cup F_1\cup F_2 \cup F_3,\tp)$: those who come from crossings, and those who come from associators. Let us note $A$ the contributions coming from associators. We will show that $A\simeq 0$.

Let us consider the reduced Kontsevich integral consisting in forgetting the contributions of the crossings. We can check that this integral is an isotopy invariant, independant of the choice of marked points because it contains only degree 2 diagrams which commute in our quotient space of diagrams.
Moreover, this reduced integral is independant of the sign of the crossings. Hence, it is an homotopy invariant. 
This proves that it reduces to $\frac{1}{48}\twowheel$ for each component, which is equivalent to 0.
We have proved that $A\simeq 0$. Then we have to take care only on the crossings.

Let us fix a component $C$ and let us follow it from its marked point. It meets an ordered sequence of components $C_1,\ldots,C_k$. Each component appear twice or never. Two cases may occur: either the component appear twice with opposite signs or it appears twice successively with the same sign.

The Campbell-Haussdorff formula tells us that each pair $\{C_i,C_j\}$ with $i<j$ creates a diagram $\YGraph$ with a coefficient $\pm\frac{1}{8}$. Moreover all $\YGraph$-shaped diagrams are obtained in that way as $A\simeq 0$.
We remark also that as $C$ does not cross itself, there are no $\multimap$-shaped diagrams. The framing of $C$ do not create $\YGraph$ or $\multimap$ diagrams. Hence, it can be ignored.

The contributions of the crossings to the $\YGraph$ diagrams are packed 4 by 4 depending on which component they link. An enumeration of different possibilities shows that the resulting coefficient is always half an integer (see figure \ref{enumer}).

\begin{figure}[htbp]
\begin{center}
\input{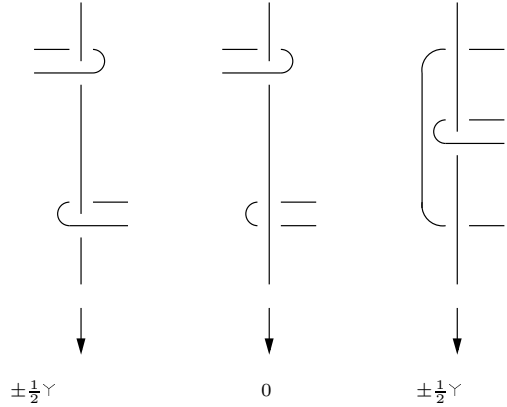}
\caption{Examples of contributions}
\label{enumer}
\end{center}
\end{figure}

\end{proof}

\begin{rem}\label{casson}
If $(M,K)$ is a pair formed of a knot in a homology sphere, then the quantity $\lambda(M)$ defined by $z_2(M,K)|_{t=1}=\frac{\lambda(M)}{2}\ThetaGraph$ is an invariant of $M$, integer if $M$ is an integer homology sphere. It is the Casson invariant of $M$.
\end{rem}

\begin{proof}
All pairs $(M,U)$ for $M$ an integer homology sphere and $U$ a trivial knot are linked by surgeries on claspers. Moreover, $z_2(M,U)$ may be written as $\frac{\lambda(M)}{2}\ThetaGraph$ for some $\lambda(M)\in\Q$. But $\lambda(S^3)=0$, hence using proposition \ref{claspsimp} and lemma \ref{demientier}, we deduce that 
$\lambda(M)\in\Z$ for all integer homology sphere $M$.

\end{proof}
\subsection{Integrality properties of $z_2$}

Let us prove now some weakened integrality property for the 2-loop part of the Kontsevich integral. This proposition is mostly important because it gives an algorithmic way to check Rozansky integrality conjecture.

\begin{prop}\label{integrality}
For all pairs $(M,K)$ formed of a knot in an integer homology sphere, the 2-loop part $z_2(M,K)$ belongs to $\frac{1}{12}\boI^{\Delta}_2$.

\end{prop}

\begin{proof}

Let us start from a surgery presentation of pair $S$-equivalent to $(M,K)$ given by a link $\boL=\{L_1,\ldots,L_n\}$ as in the second part of lemma \ref{projection}.
We note that each component bounds a disc naturally decomposed in an upper hemisphere and a lower hemisphere. We choose a base point and an arc from the base point to each component. We choose a marked point at the left side of the equatorial line of each disc.

As we are interested in the 2-loop part, we will evaluate only diagrams which appear on the figure \ref{white7}. In particular, any connected diagram containing at least 3 trivalent vertices will be considered equal to 0.

\mbox{\bf Contribution of the framing}

Ley $L_i'$ be the component $L_i$ without twists. We can compute $\sigma Z(\hat{L}_i,\tp)$ from $\sigma Z(\hat{L'}_i,\tp)$ thanks to a wheeling map using a method already used in \cite{kr2}.

Let $W$ be the degree 1 part of $\sigma(\hat{L'}_i,\tp)$. It is integer since $L'i$ has no self intersection and each pair of discs meet in an even number of points. Let $w_i$ be the framing of the $i$-th component.

We recall that there is a normalization factor $\exp(\frac{1}{48}\twowheel)$ for each component symbolized by the notation $\hat{L}$. From  $L'$ to $L$, we must multiply the $i$-th component of  $Z(\hat{L'}_i,\tp)$ with $\exp_{\#}(\frac{w_i}{2}\isolatedchord_i)$ for all $i$.

As this element contains legs colored by $i$, we have (see \cite{kr2} or \cite{th})

$$\sigma Z(L_i,\tp)=\prod_i\partial_{\Omega_i}
\left(\prod_i\partial^{-1}_{\Omega_i}
\sigma Z(L'_i,\tp) 
\exp(
\sum_i\frac{w_i}{2}\alg{i}{\arete}{i}+
\frac{1}{48}\sum_i\alg{i}{\twowheel}{i}
-\frac{1}{24}\sum_i w_i\ThetaGraph)
\right).$$

This last quantity is equal to
$$\sigma Z(L'_i,\tp)\exp\left(\sum_i (\frac{w_i}{2}\alg{i}{\arete}{i} + \frac{1}{24}\sum_{j\ne i}w_i\,\,\overset{w_{i,j}}{\alg{i}{\twowheel}{j}}+\frac{w_i^2+1}{48}\alg{i}{\twowheel}{i})\right).$$
This allows us to replace the computation of $\sigma Z(\hat{L}_i,\tp)$ by the computation of $\sigma Z(L'_i,\tp)$.

To compute the integral, we use heavily the fact that $L'_i$ is an element of the class $\boC$. It allows us:
\begin{itemize}
\item
to eliminate contributions coming from diagrams $A$ and $B$ of figure \ref{classec}.
\item
to separate contributions linking 1, 2 and more components.
\end{itemize}

Let $\boH$ be the abelian group generated by the diagrams $\frac{1}{12}\multimap\multimap$, $\frac{1}{12}\multimap\YGraph$, $\frac{1}{12}\YGraph\YGraph$, $\frac{1}{12}\twowheel$, $\frac{1}{12}I$ et $\frac{1}{12}Y_2$ colored by monomials and with legs parametrized by the set $\{1,\ldots,n\}$.

It is sufficient to compute $\sigma Z(L'_i,\tp) \mod \boH$ because we know that diagrams of $\boH$ produce a denominator $\frac{1}{12}$ and hence satisfy the theorem.

\mbox{\bf Elimination of double crossing and associators:}

It is clear that the constribution of the tangle $A$ of figure \ref{classec} is $\frac{1}{2}Y_2$ where $Y_2$ is a diagram linking the four component represented in the figure. This diagram lies in $\boH$.

The same computation holds for the associator, except that the contribution is a $\YGraph$-graph whose legs are attached to an edge. There are 4 internal trivalent vertices in this diagram, so that it may be ignored.

\mbox{\bf Study of the elementary link:}

Each isolated component is trivial, which proves that there is only one diagram colored by a single component $L_i$: $\frac{1}{48}\alg{i}{\twowheel}{i}$.

Each pair of linked components form what we call an elementary link. To evaluate the contributions linked to two components, we compute explicitely the Kontsevich integral of the tangle presented on the figure \ref{maillon}.

\begin{figure}[htbp]
\begin{center}
\input{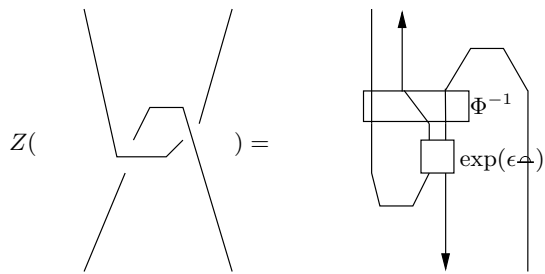}
\caption{Elementary link}
\label{maillon}
\end{center}
\end{figure}
Let $\epsilon$ be the sign of the linking of the two components oriented as in the figure. We note $a$ and $b$ the two components. A direct computation shows that the part of $Z(L_a\cup L_b,\tp)$ linked to both components is precisely:

$$\exp(\epsilon \alg{a}{\arete}{b}+\frac{\epsilon^2}{8}\alg{a}{\twowheel}{b}+\frac{\epsilon+2\epsilon^3}{24}H^a_b).$$
The notation $H^a_b$ designs the graph $Y_2$ whose left legs are colored by $a$ and right legs are colored by $b$.

\mbox{\bf Contributions linked to more than 2 components}

As the tangles $A$ and $B$ disappear, it is sufficient to know how to combine the integrals of the elementary links to obtain the integral of $\boL$.
For that, we can ignore the associators necessary to glue the elementary links together. It only remains to multiply all these elements along the components. We will perform this computation by using a generalized Baker-Campbell-Haussdorff formula.

Let $H(a_1,\ldots,a_p)$ be the map from $\boB(a_1,\ldots,a_p,S)$ to $\boB(a,S)$ which is obtained by gluing on $a_1,\ldots,a_p$ the following graph:

\begin{eqnarray*}
H(a_1,\ldots,a_p)=\exp(\sum_i a_i+\frac{1}{2}\sum_{i<j}[a_i,a_j]+\frac{1}{6}\sum_{i<j<k}[a_i,[a_j,a_k]]+\frac{1}{6}\sum_{i<j<k}[[a_i,a_j],a_k]\\
+\frac{1}{12}\sum_{i<j}[a_i,[a_i,a_j]]+\frac{1}{12}\sum_{i<j}[a_j,[a_j,a_i]]).
\end{eqnarray*}

This operator describes how to glue all elementary links along a same component. We distinguish two contributions of the elementary links: connected diagrams with 2 internal vertices which cannot be glued to $H$ and edges which have an integer coefficient.

The diagrams of $H$ having two internal vertices already have a denominator dividing 12. After integration, they will produce elements of $\boH$.

The contributions obtained by gluing two distinct $\YGraph$ has denominator $\frac{1}{4}$ except for the terms $\frac{\epsilon^2}{8}\alg{a}{\twowheel}{b}$ which are already counted in the elementary link.

Nevertheless, each graph $\frac{1}{2}\YGraph$ glued to itself produce a denominator $\frac{1}{8}$ which has to be computed more precisely.

Let $w_{i,j}$ be the equivariant linking number of the $i$-th and the $j$-th component. By construction, it is either 0 either an invertible element of $\Z[t,t^{-1}]$. We will note  $|\epsilon t^k|=|\epsilon|t^k$. 

We sum up the preceding facts in the following formula:

\begin{align*}\label{zlprime}
\sigma Z(L',\tp)=\exp\left(\sum_{\{i,j\}} w_{ij}\alg{i}{\arete}{j}+\sum_i \frac{1}{48}\alg{i}{\twowheel}{i}+\sum_{\{i,j\}}\frac{|w_{ij}|}{8}\alg{i}{\twowheel}{j}+\sum_{\{i,j\}}\frac{w_{ij}}{8}H^i_j+\sum_{\{i,j,k\}}\frac{1}{2}\vphantom{\YGraph_i}^{w_{ij}}\YGraph_i^{w_{ik}}\right)\\ \mod \boH
\end{align*}

\mbox{\bf Conclusion}

To take care of the framing of the components and of their normalization factors, we have to multiply 
$\sigma Z(L'_i,\tp) \mod \boH$ by the quantity:
$$\exp\left(\sum_i (\frac{w_i}{2}\alg{i}{\arete}{i} + \frac{w_i}{24}\sum_{j\ne i}\overset{w_{i,j}}{\alg{i}{\twowheel}{j}}+\frac{w_i^2+1}{48}\alg{i}{\twowheel}{i})+\frac{\sigma}{16}\ThetaGraph\right)$$
In this formula, $\sigma$ is the signature of the non-equivariant linking matrix of $L$.

To prove the integrality statement about the 2-loop part, we will integrate the preceding formula and compute it modulo diagrams such as $\frac{1}{12}\theta(\frac{P}{\Delta},\frac{Q}{\Delta},\frac{R}{\Delta})$ with $P,Q,R\in\Z[t,t^{-1}]$.

We also know that the 2-loop part for $t=1$ is half-integer (see remark \ref{casson}). It means that we may also quotient by diagrams such as $\alpha \ThetaGraph$ pour $\alpha \in\Q$. We note $\boH'$ the subgroup of $\boA^{\Delta}_2$ generated by these diagrams and $\frac{1}{12}\boI^{\Delta}_2$

Then, if we prove that $z_2(M,K)=x+\alpha\ThetaGraph$ for $x\in\frac{1}{12}\boI^{\Delta}_2$, it suffices to take $t=1$ to prove that $\alpha\in\frac{1}{12}\Z$ and hence $z_2(M,K)\in\frac{1}{12}\boI^{\Delta}_2$.

Let us compute each contribution $\mod \boH'$ separately.
\begin{itemize}\item
Case of $\sum_i\frac{w_i^2+2}{48}\alg{i}{\twowheel}{i}+\frac{w_i}{24}\sum_{j\ne i}\overset{w_{i,j}}{\alg{i}{\twowheel}{j}}+\frac{\sigma}{16}\ThetaGraph$.

Modulo $\boH'$, we can ignore the term $\frac{\sigma}{16}\ThetaGraph$.
Let $D$ be the diagonal matrix whose coefficients are the $w_i$.
This part gives the contribution: 
$\theta(-\frac{1}{48}\Tr (D^2+2\id) W^{-1}- \frac{1}{24}\Tr D (W-D) W^{-1})=\theta(-\frac{1}{24}\Tr(W+W^{-1})+\frac{1}{48}\Tr D^2W^{-1})$.

But $\Tr(W+W^{-1})=\frac{Q}{\Delta}$ where $Q$ is a symetric polynomial. The non constant monomials come into pairs and vanishes modulo $\boH'$. The only remaining diagrams are the integer multiples of the diagram $\frac{1}{24}\theta(\frac{1}{\Delta})$ which also live in $\boH'$.
The contribution $\frac{1}{48}\theta(\Tr D^2W^{-1})$ does not leave in $\boH'$.

\item
Case of $\sum_{\{i,j\}}\frac{|w_{ij}|}{8}\alg{i}{\twowheel}{j}+\sum_{\{i,j\}}\frac{w_{ij}}{8}H^i_j$.

We integrate these diagrams modulo $\boH'$.

\begin{align*}
\irat \left[\frac{|w_{ij}|}{8}\alg{i}{\twowheel}{j}+\frac{w_{ij}}{8}H^i_j\right] \\
=\frac{1}{8}\left[
\theta(w_{ij}w_{ji}^{-1})+\alg{w_{ij} w_{ji}^{-1}}{\dumbbell}{w_{ij}w_{ji}^{-1}}
-\theta(w_{ij}w_{ji}^{-1},w_{ji}w_{ij}^{-1})
+w_{ij}(1)\theta(w_{ii}^{-1},w_{jj}^{-1}) \right]\\
=\frac{1}{8}\left[\theta (w_{ij}w_{ji}^{-1}) + \theta(w_{ij}w_{ji}^{-1},w_{ij}w_{ji}^{-1})+w_{ij}(1)w_{ii}^{-1}(1)w_{jj}^{-1}(1)\ThetaGraph\right]\\
=0\mod \boH'
\end{align*}
This last equality is a consequence of the formula  $\Delta U = U\otimes U \mod 2$ for all polynomial $U\in\Z[t,t^{-1}]$.

\item
Case of $\sum_{\{i,j,k\}}\frac{1}{2}(\frac{1}{2}\vphantom{\YGraph_i}^{w_{ij}}\YGraph_i^{w_{ik}})^2$.

\begin{lem}
Let $\YGraph$ be a graph colored by three elements, $a_i$, $a_j$ and $a_k$ of $\Z[t,t^{-1}]$. These elements lie on the three edges in positive order, and are directed through the exterior of the graph.

Let $W$ be a matrix of $\Z[t,t^{-1}]$ whose determinant is equal to $\pm 1$ for $t=1$ and note $\epsilon P=P(1)$ for $P\in \Z[t,t^{-1}]$. Then, modulo 2, $\irat \YGraph^2= \epsilon(w^{-1}_{ii}a_i\ba{a_i})\theta(a_jw^{-1}_{jk}\ba{a_k})+\epsilon(w^{-1}_{jj}a_j\ba{a_j})\theta(a_kw^{-1}_{ki}\ba{a_i})+\epsilon(w^{-1}_{kk}a_k\ba{a_k})\theta(a_iw^{-1}_{ij}\ba{a_j})+\alpha\ThetaGraph$ for some $\alpha\in\Z/2\Z$.

\end{lem}
The proof of this lemma consist in enumerating the 15 ways for gluing two copies of  $\YGraph$. We eliminate 9 of them for symetry reasons, and we simplify the last 6 thanks to the (IHX) relations and the sliding relations modulo 2.

We apply this lemma directly to find for all triple $\{i,j,k\}$:

$\irat(\vphantom{\YGraph_i}^{w_{ij}}\YGraph_i^{w_{ik}})^2 = \epsilon(w^{-1}_{ii})\theta(w_{ij}w^{-1}_{jk}w_{ki})+
\epsilon(w^{-1}_{jj}w_{ij}w_{ji})\theta(w^{-1}_{ki}w_{ik})+
\epsilon(w^{-1}_{kk}w_{ik}w_{ki})\theta(w^{-1}_{ji}w_{ij}) +\alpha\ThetaGraph\mod 2$.
The two last terms are linked by the permutation $(jk)$.

Each term may be written as $\phi(i,j,k)$ and we compute the sum 
\begin{multline*}
\kern2em\sum_{\{i,j,k\}}\phi(i,j,k)=
\sum_i \frac{1}{2}\sum_{j,k,j\ne i,k\ne i,j\ne i}\phi(i,j,k)\\
=\sum_i \frac{1}{2}\left[\sum_{j,k}\phi(i,j,k)-\sum_{j}\phi(i,j,i)-\sum_k \phi(i,i,k)-\sum_l \phi(i,l,l)+2\phi(i,i,i)\right]
\end{multline*}

If $\phi(i,j,k)=\epsilon(w^{-1}_{ii})\theta(w_{ij}w^{-1}_{jk}w_{ki})$, then $\sum_{j,k}\phi(i,j,k)=\sum_j \phi(i,j,i)=\sum_k\phi(i,i,k)=\epsilon(w^{-1}_{ii})w_{ii}$.

Hence, $\sum_l\phi(i,l,l)=\sum_l\epsilon(w^{-1}_{ii})\theta(w_{il}w_{ll}^{-1}w_{li})$ and $\phi(i,i,i)=\epsilon(w^{-1}_{ii})\theta(w_{ii}^2w^{-1}_{ii})$.

The contribution of this term modulo $\boH'$ is $-\frac{1}{16}\sum_{i,l}\epsilon(w^{-1}_{ii})\theta(w_{il}w_{ll}^{-1}w_{li})$.

If $\phi(i,j,k)=\epsilon(w^{-1}_{jj}w_{ij})\theta(w^{-1}_{ki}w_{ik})$ then the two remaining contributions may be written as:

$\sum_i \left[\sum_{j,k}\phi(i,j,k)-\sum_{j}\phi(i,j,i)-\sum_k \phi(i,i,k)-\sum_l \phi(i,l,l)+2\phi(i,i,i)\right]$.

But
\begin{align*}
\sum_{j,k}\phi(i,j,k)=\sum_j\epsilon(w^{-1}_{jj}w_{ij})\text{ and }
\sum_{k}\phi(i,i,k)=\epsilon(w^{-1}_{ii}w_{ii})\\
\sum_j\phi(i,j,i)=\sum_j\epsilon(w_{jj}^{-1}w_{ij})\theta(w^{-1}_{ii}w_{ii})\text{ and }
\sum_l\phi(i,l,l)=\sum_l\epsilon(w_{ll}^{-1}w_{il})\theta(w^{-1}_{li}w_{il})
\end{align*}
But $w_{il}$ is a  monomial so that we have $w_{il}=\epsilon(w_{il})w_{il} \mod 2$.

Hence, $\sum_{i,l}\phi(i,l,l)=\sum_{i,l}\epsilon(w_{ll}^{-1})\theta(w^{-1}_{li}w_{il})=\sum_l\epsilon(w^{-1}_{ll})\in \boH'$.
\end{itemize}

The remaining terms are the following:
\begin{eqnarray*}
-\frac{1}{16}\sum_{i,l}\epsilon(w^{-1}_{ii})\theta(w_{il}w_{ll}^{-1}w_{li})
+\frac{1}{48}\sum\theta(w_{ii}^2w^{-1}_{ii}) \\
=\frac{1}{48}\sum_{i,l}\epsilon(w^{-1}_{ii})\theta(w_{il}w_{ll}^{-1}w_{li})
+\frac{1}{48}\sum\theta(w_{ii}^2w^{-1}_{ii}) \mod \boH'
\end{eqnarray*}

But we have $w_{il}w_{li}=\epsilon(w_{il}w_{li})$, so we can replace it to obtain

$\sum_i \epsilon(w_{li})\epsilon(w_{ii}^{-1})\epsilon(w_{il})=\sum_{i,j}\epsilon(w_{li})\epsilon(w_{ij}^{-1})\epsilon(w_{jl})= w_{ll}\mod 2$.
Hence, as we have $w_{ii}^2=w_{ii}\mod 2$, the two terms vanish, and allows to prove that $\irat\sigma Z(\hat{L}_i,\tp)\in \boH'$. This prove finally the integrality result about 2-loop part.

\end{proof}

\section{Remarks and applications}



The proposition \ref{integrality} can be stated in the following more refined form: 
\begin{prop}
Let $n$ be a positive integer. We call monomial a $n\times n$ hermitian matrix with monomial coefficients in  $\Z[t,t^{-1}]$ invertible over $\Z$. Let $\boM_n$ be the set of such matrices. There is an algebraic map  $\phi_n:\boM_n\to\boA_2$ such that for all $W\in \boM_n$ and all pairs $(M,K)$ $S$-equivalent to $W$, we have $$z_2(M,K)=\phi_n(W) \mod \frac{1}{2}\boI^{\Delta}_2.$$
\end{prop}

The formula for $\phi_n$ is huge and not enlightening: it contains contributions of elementary links, self-crossings and of the generalized Campbell-Hausdorff formula. We programmed it to check the validity of Rozansky's conjecture on some monomial matrices.

\begin{quest}
Find a simple formula for $\phi_n$, use it to check, prove or disprove Rozansky's conjecture.
\end{quest}

Let us prove this proposition communicated to us by Y. Tsustsumi as an application of our results.
\begin{prop}
Let $(M,K)$ be a pair formed of a knot $K$ in an integer homology sphere $M$, and $r$ an integer such that $\Sigma^r(M,K)$ is again a knot in an integer homology sphere. 
Then for all pair $(M',K')$ $S$-equivalent to $(M,K)$, $r$ divides $\lambda(\Sigma^r(M,K))-\lambda(\Sigma^r(M',K'))$.

If Rozansky's conjecture is true, we have $\lambda(\Sigma^r(M,K))\equiv \frac{1}{8}\sigma_r(M,K)\mod r$.
\end{prop}
\begin{proof}
Let $\Delta$ be the normalized Alexander polynomial of $(M,K)$ and $(M',K')$. 
Thanks to lemma \ref{demientier}, $\frac{1}{2}x=z_2(M',K')-z_2(M,K)\in\frac{1}{2}\boI^{\Delta}_2$.
We use the formula of \cite{lift} to obtain $\frac{1}{2}\lambda(\Sigma^r(M,K))=\frac{\sigma_r(M,K)}{16}+\lift_r z_2(M,K)|_{t=1}$ and $\frac{1}{2}\lambda(\Sigma^r(M',K'))=\frac{\sigma_r(M',K')}{16}+\lift_r z_2(M,K')|_{t=1}$.

As $(M,K)$ and $(M',K')$ are $S$-equivalent, we have $\sigma_r(M,K)=\sigma_r(M',K')$. To prove the assertion, we only have to show that $r$ divides $\lift_r x|_{t=1}$.
In order to compute $\lift_r x$, we only need to find some polynomial $P\in\Z[t,t^{-1}]$ such that $\frac{1}{\Delta(t)}=\frac{P(t)}{\Delta_r(t^r)}$, where $\Delta_r$ is the Alexander polynomial of $\Sigma^r(M,K)$. The result will follow from the definition of the $\lift_r$ map.

But $\Delta_r$ is a normalized polynomial with integer coefficients because the underlying manifold of the pair $\Sigma^r(M,K)$ is an integer homology sphere. Its roots are the $r$-th powers of the roots of $\Delta$. Hence, $\Delta$ divides $\Delta_r$ over $\Q$.

For any polynomial $P$, we note $c(P)$ the content of $P$, i.e. the greatest common divisor of its coefficient. As $\Delta(1)=\Delta_r(1)=1$, we have $c(\Delta)=c(\Delta_r)=1$. If we write $P=\frac{a}{b}Q$ for some irreducible fraction $\frac{a}{b}$ and $Q$ an integer polynomial with $c(Q)=1$, we deduce from $\Delta_r(t^r)=P(t)\Delta(t)$ that $b\Delta_r(t^r)=a Q(t)\Delta(t)$. By multiplicativity of the contents, we find $b=\pm a$ which proves that $P$ is an integer polynomial.
\end{proof}

\bibliographystyle{halpha}
\bibliography{integ}

\end{document}